\documentclass[12pt]{amsart}
\usepackage{times}
\usepackage{amsfonts}
\usepackage{amsthm}
\usepackage{amsmath}
\advance \topmargin by -\headheight
\advance \topmargin by -\headsep
\evensidemargin \oddsidemargin
\newtheorem{theorem}{Theorem}[section]
\newtheorem{lemma}[theorem]{Lemma}
\newtheorem{corollary}[theorem]{Corollary}

\newtheorem{remark}[theorem]{Remark}

\newtheorem{definition}[theorem]{Definition}

\begin{document}

\title{A Free Entropy Dimension Lemma}

\author{Kenley Jung}

\address{Department of Mathematics, University of California,
Berkeley, CA 94720-3840,USA}

\email{factor@math.berkeley.edu}
\thanks{Research supported by the NSF Graduate Fellowship Program}

\begin{abstract} Suppose $M$ is a von Neumann algebra with normal, tracial
state $\varphi$ and $\{a_1, \ldots,a_n\}$ is a set of self-adjoint
elements in $M.$ We provide an alternative uniform packing description of
$\delta_0(a_1,\ldots,a_n),$ the modified free entropy dimension of
$\{a_1,\ldots,a_n\}.$ \end{abstract}

\maketitle

\begin{center}
\emph{For Arlan Ramsay}
\end{center}

\vspace{.1in}

	In the attempt to understand the free group factors Voiculescu
created a type of noncommutative probability theory.  One facet of the
theory involves free entropy and free entropy dimension, applications of
which have answered some old operator algebra questions ([1], [4]).  
Roughly speaking, given self-adjoint elements $a_1,\ldots,a_n$ in a von
Neumann algebra $M$ with normal, tracial state $\varphi$ a matricial
microstate for $\{a_1,\ldots,a_n\}$ is an $n$-tuple of self-adjoint $k
\times k$ matrices which together with the normalized trace, approximate the algebraic behavior of the $a_i$ under
$\varphi.$ Taking a normalization of the logarithmic volume of such
microstate sets followed by multiple limiting processes yields a number,
$\chi(a_1,\ldots, a_n),$ called the free entropy of $\{a_1,\ldots,a_n\}.$
One can think of free entropy as the logarithmic volume of the $n$-tuple.  
The (modified) free entropy dimension of $\{a_1,\ldots,a_n\}$ is

\[ \delta_0(a_1,\ldots,a_n) = n+ \limsup_{\epsilon \rightarrow 0} 
\frac{\chi(a_1+\epsilon s_1,\ldots,a_n+\epsilon s_n:s_1,\ldots,s_n)}{|\log 
\epsilon|}\]

\noindent where $\{s_1,\ldots,s_n\}$ is a semicircular family freely
independent with respect to $\{a_1,\ldots, a_n\}$ and $\chi(:)$ is a 
technical modification of $\chi$ (see [4]). 

\vspace{.1in}

	Free entropy dimension was inspired by Minkowski dimension.  
Recall that for a subset $A \subset \mathbb R^d$ the (upper) Minkowski
dimension of $A$ is \[ d + \limsup_{\epsilon \rightarrow 0} \frac {\log
\lambda (\mathcal N_{\epsilon}(A))}{|\log \epsilon|}\] where $\lambda$
above denotes Lebesgue measure and $\mathcal N_{\epsilon}(A)$ is the
$\epsilon$-neighborhood of $A.$ Minkowski dimension has an equivalent
formulation in terms of uniform packing dimension.  The (upper) uniform
packing dimension of $A$ is \[ \limsup_{\epsilon \rightarrow 0} \frac{
\log P_{\epsilon}(A)}{|\log \epsilon|} = \limsup_{\epsilon \rightarrow 0}
\frac {\log K_{\epsilon}(A)}{|\log \epsilon|} \] where $A$ is endowed with
the Euclidean metric, $P_{\epsilon}(A)$ is the maximum number of elements in a 
collection of mutually disjoint open $\epsilon$ balls of $A,$ and
$K_{\epsilon}(A)$ is the minimum number of open $\epsilon$ balls required
to cover $A$ (the quantities above make sense in the setting of an
arbitrary metric space).  It is easy to see that the Minkowski dimension
and the uniform packing dimension of $A$ are always equal.

\vspace{.1in}

	In this paper we present a lemma which formulates a similar
metric description of $\delta_0$: free entropy dimension can be
described in terms of packing numbers with balls of equal radius.

\vspace{.1in}

	The alternative description comes as no surprise in view of 
both the definition of $\delta_0$ and the techniques in 
estimations thereof.  The proof follows the
classical one with the addition of the properties of $\chi$ proven in 
[3] and the strengthened asymptotic freeness results of [5].

\section{Preliminaries} Throughout $M$ is a von Neumann algebra with
normal, tracial state $\varphi$ and $\{a_1,\ldots,a_n\}$ is a set of
self-adjoint elements in $M.$ We use the symbols $\chi$ and $\delta_0$ to
designate the same quantities introduced in [4].  $M^{sa}_k(\mathbb C)$
denotes the set of $k \times k$ self-adjoint complex matrices and
$(M^{sa}_k(\mathbb C))^n$ is the set of $n$-tuples with entries in
$M^{sa}_k(\mathbb C).$ $tr_k$ is the normalized trace on the $k \times k$
complex matrices.  $\|\cdot \|_2$ is the inner product norm on
$(M^{sa}_k(\mathbb C))^n$ given by the formula $\|(x_1,\ldots,x_n)  
\|_2^2 = \sum_{i=1}^n k \cdot tr_k(x_i^2)$ and vol denotes Lebesgue
measure with respect to the $\| \cdot \|_2$ norm.  For any $k \in \mathbb
N$ denote by $\rho_k$ the metric on $(M^{sa}_k(\mathbb C))^n$ induced by
the norm $ k^{ - \frac{1}{2}} \cdot\| \cdot \|_2.$ For a metric space
$(X,d)$ and $\epsilon >0$ write $P_{\epsilon} (X, d)$ for the maximum
number of elements in a collection of mutually disjoint open $\epsilon$
balls of $X$ and $K_{\epsilon} (X, d)$ for the minimum number of open
$\epsilon$ balls required to cover $X.$ Observe that $P_{\epsilon}(X,d)  
\geq K_{2 \epsilon}(X, d) \geq P_{4 \epsilon}(X, d).$ Finally for $S
\subset X$ denote by $\mathcal N_{\epsilon}(S)$ the
$\epsilon$-neighborhood of $S.$

\section{The Lemma}

\begin{definition}  For any $k, m \in \mathbb N,$ and $R, \gamma, \epsilon 
>0$ define successively 

\[ \mathbb P_{\epsilon, R}(a_1, \ldots, a_n;m,k,\gamma) = 
P_{\epsilon}(\Gamma_R(a_1, \ldots,a_n; m,k,\gamma) , \rho_k),\]  

\[ \mathbb P_{\epsilon, R}(a_1, \ldots, a_n;m,\gamma) = \limsup_{k 
\rightarrow \infty} k^{-2} \cdot \log (\mathbb P_{\epsilon, R}(a_1, 
\ldots,a_n ; m, k, 
\gamma)), \]

\[ \mathbb P_{\epsilon, R}(a_1, \ldots, a_n) = \inf \{ \mathbb 
P_{\epsilon, 
R}(a_1, \ldots,a_n; m, \gamma) : m \in \mathbb N, \gamma >0 \}, \]

\[ \mathbb P_{\epsilon}(a_1,\ldots,a_n) = \sup_{R > 0} \{ \mathbb 
P_{\epsilon, R}(a_1, \ldots, a_n) \}.\]

\end{definition}

\noindent{\it Remark.} If $b_1, \ldots, b_p \in M,$ then 
define $\mathbb P_{\epsilon}(a_1,\ldots, a_n: b_1, \ldots, b_p)$ to be the quantity obtained by replacing $\Gamma_R(a_1,\ldots,a_n;m,k,\gamma)$ in the definition with $\Gamma_R(a_1,\ldots,a_n:b_1,\ldots,b_p;m,k,\gamma).$  Similarly, we define 
$\mathbb K_{\epsilon}(a_1, \ldots, a_n)$ and all its associated quantities 
by replacing $P_{\epsilon}$ in the first line of Definition 2.1 with 
$K_{\epsilon}.$  Define $\mathbb 
K_{\epsilon}(a_1,\ldots,a_n:b_1,\ldots,b_p)$ in the same way  $\mathbb 
P_{\epsilon}(a_1,\ldots,a_n:b_1,\ldots,b_p)$ was defined.

\vspace{.1in} 

For any self-adjoint elements $h_1,\ldots,h_n \in M$ denote
by $\underline{\chi} (h_1,\ldots,h_n)$ the number obtained by replacing
the $\limsup$ in the definition of $\chi$ with $\liminf.$ $\mathbb
P_{\epsilon}(\cdot)$ being a normalized limiting process of the
logarithmic microstate space packing numbers, we have:

\begin{lemma} If $\{h_1,\ldots,h_n \}$ is a set of self-adjoint elements
in $M$ which is freely independent with respect to $\{a_1,\ldots,a_n\}$
and $\underline{\chi}(h_1,\ldots,h_n) > - \infty,$ then

\begin{eqnarray*} n + \limsup_{\epsilon \rightarrow 0} \frac {\chi(a_1 +
\epsilon h_1, \ldots, a_n + \epsilon h_n:h_1,\ldots, h_n)}{|\log
\epsilon|} & = & \limsup_{\epsilon \rightarrow 0} \frac {\mathbb
K_{\epsilon}(a_1, \ldots, a_n)}{| \log \epsilon|} \\ & = &
\limsup_{\epsilon \rightarrow 0} \frac {\mathbb P_{\epsilon}(a_1, \ldots,
a_n)}{| \log \epsilon|}.  \end{eqnarray*} \end{lemma} 

\begin{proof} Clearly it suffices to show equality of the first and last
expressions above since $P_{\epsilon}(\cdot) \geq K_{2 \epsilon}(\cdot)  
\geq P_{4 \epsilon}(\cdot).$ Furthermore, we can assume that
$\{a_1,\ldots,a_n\}$ has finite dimensional approximants since the
equalities hold trivially otherwise.  Set $C = \max \{ \|h_i\| \}_{1 \leq
i \leq n}.$ First we show that the free entropy expression is greater
than or equal to the $\mathbb P_{\epsilon}$ expression.  Suppose $0 <
\epsilon < \frac{1}{2(C+1)},$ $m \in \mathbb N$ with $m >n,$ $\gamma >
0,$ and $R > \max\{\|a_i\|\}_{1 \leq i \leq n}.$ 

\vspace{.1in}

Corollary 2.14 of [5] provides an $N\in \mathbb N$ such that if $k \geq N$
and $\sigma$ is a Radon probability measure on $((M_{k}^{sa}(\mathbb
C))_{R+1})^{2n}$ (the subset of $(M^{sa}_k(\mathbb C))^{2n}$ consisting of
$2n$-tuples whose entries have operator norm no greater than $R+1$)  
invariant under the $U_k$-action

\[(\xi_1,\ldots,\xi_n,\eta_1,\ldots,\eta_n) \mapsto
(\xi_1,\ldots,\xi_n,u\eta_1 u^{*},\ldots,u\eta_n u^{*}),\] 

\noindent then $\sigma(\omega_k) > \frac {1}{2}$ where $\omega_k$ is

\[ \{ (\xi_1,\ldots,\xi_n,\eta_1,\ldots,\eta_n)\in((M_k^{sa}(\mathbb
C))_{R+1})^{2n} \hspace{-.05in} :  \hspace{-.05in}\{\xi_1,\ldots,\xi_n\}
\hspace{.05in} \text{and} \hspace{.05in} \{\eta_1,\ldots,\eta_n\}
\text{are} \left (m, \gamma/4^m \right) \text{ - free} \}.  \]

\noindent With respect to the $\rho_k$ metric for each $k$ find a
collection of mutually disjoint open $ 2 C \epsilon \sqrt{n}$
balls of $\Gamma_R (a_1, \ldots, a_n; m, k, 
\gamma/(8(R+2))^m )$ of maximum cardinality and denote the
centers of these balls by $<(x_{1j}^{(k)}, \ldots, x_{nj}^{(k)})>_{j \in
S_k}$.  Define $\mu_k$ to be the uniform atomic probability measure
supported on $<(x_{1j}^{(k)}, \ldots, x_{nj}^{(k)})>_{j \in S_k}$ and
$\nu_k$ to be the probability measure obtained by restricting vol to
$\Gamma_{2 C \epsilon} (\epsilon h_1, \ldots, \epsilon h_n; m,k, \gamma/8^m)$ and normalizing appropriately.  $\mu_k \times
\nu_k$ is a Radon probability measure on $((M_k^{sa}(\mathbb
C))_{R+1})^{2n}$ invariant under the $U_k$-action described above.  So for
$k \geq N$ $(\mu_k \times \nu_k)(\omega_k) > \frac{1}{2}.$

\vspace{.15in}

For $k\in \mathbb N$ and $j \in S_k$ define $F_{jk}$ to be the set of 
all $(y_1, 
\ldots, y_n) \in \Gamma_{2C \epsilon} (\epsilon h_1, \ldots, \epsilon 
h_n ; 
m,k, \gamma/8^m )$ such that $(y_1, \ldots, y_n)$ and 
$(x_{1j}^{(k)}, \ldots, x_{nj}^{(k)})$ are $\left (m, \frac {\gamma}{4^m} 
\right)$-free.  

\[ \frac{1}{2} < (\mu_k \times \nu_k)(\omega_k) = 
\sum_{j \in S_k} \frac {1}{|S_k|} \cdot \nu_k(F_{jk}) = 
\sum_{j \in S_k} \frac {1}{|S_k|} \cdot \frac{\text {vol} (F_{jk})}{\text{vol} 
\left(\Gamma_{2C\epsilon } (\epsilon h_1,\ldots, 
\epsilon h_n ; m,k, \gamma/8^m )\right)}. \]

\noindent It follows that for $k \geq N$
\[ \frac{1}{2} \cdot |S_k| \cdot \text{vol} 
(\Gamma_{1} (\epsilon h_1, \ldots, 
\epsilon h_n; m, k, \gamma /8^m ) ) 
< \sum_{j \in S_k} \text{vol}(F_{jk}). \]

\noindent Set $E_{jk} = (x_{1j}^{(k)}, \ldots, x_{nj}^{(k)}) + F_{jk}.$
$F_{jk}$ is a set contained in the open ball of $\rho_k$ radius
 $2 C \epsilon \sqrt{n}$ centered at $(0, \ldots, 0).$ Thus $<E_{jk}>_{j
\in S_k}$ is a collection of mutually disjoint sets.  So \[\bigsqcup_{j
\in S_k} E_{jk} \subset \Gamma_{R+1}(a_1 + \epsilon h_1, \ldots, a_n +
\epsilon h_n : \epsilon h_1, \ldots, \epsilon h_n ; m,k, \gamma). \]

\noindent Thus, for any $\frac{1}{2(C+1)} > \epsilon >0,$ $m \in \mathbb 
N$ sufficiently large,
and $\gamma >0$ 

\begin{eqnarray*} & & \chi_{R+1}(a_1 +\epsilon h_1, \ldots, a_n+ \epsilon 
h_n:  
\epsilon h_1, \ldots, \epsilon h_n;m, \gamma) \\
& \geq & \limsup_{k \rightarrow \infty} 
\left( k^{-2} \cdot \log (\text{vol}(\bigsqcup_{j \in S_k} E_{jk})) + 
\frac {n}{2} \log k \right) \\ & \geq & 
\limsup_{k \rightarrow \infty}\left [ k^{-2} \cdot \log \left 
(\frac{1}{2} \cdot |S_k| \cdot \text{vol} \left (\Gamma_{2 C \epsilon}
(\epsilon h_1, \ldots, \epsilon h_n;m,k, \gamma / 8^m)\right )\right ) 
+ \frac {n}{2} \log k \right ] \\  & \geq & 
\limsup_{k\rightarrow \infty} \left [ k^{-2} \cdot \log (|S_k|) \right] 
+ \liminf_{k \rightarrow \infty} \left[ k^{-2} \cdot \log 
( \text{vol}(\Gamma_{2 C \epsilon}(\epsilon h_1, \ldots, \epsilon h_n;m,k,\gamma 
/ 8^m))) 
+ \frac{n}{2} \log k\right ] \\  
& \geq & \mathbb P_{2 C \epsilon \sqrt{n}, R+1} \left(a_1, \ldots, 
a_n;m, 
\gamma/ (8(R+2))^m \right) + \underline{\chi_{2 C \epsilon}}(\epsilon h_1, 
\ldots, 
\epsilon h_n) \\ 
& \geq & \mathbb P_{2 C \epsilon \sqrt{n}, R+1} (a_1, 
\ldots, a_n) + n \log 
\epsilon + \underline{\chi}(h_1,\ldots,h_n). 
\end{eqnarray*} 

\vspace{.15in}

By the chain of inequalities of the preceding paragraph it follows that
\begin{eqnarray} \nonumber \chi(a_1 + \epsilon h_1, \ldots, a_n + \epsilon
h_n: h_1, \ldots, h_n) & = & \chi(a_1 + \epsilon h_1, \ldots, a_n +
\epsilon h_n: \epsilon h_1, \ldots, \epsilon h_n) \\ \nonumber & \geq &
\mathbb P_{2 C \epsilon \sqrt{n} , R+1}(a_1, \ldots, a_n) + n \log 
\epsilon +
\underline{\chi}(h_1,\ldots,h_n).  \nonumber \end{eqnarray} This being
true for any $R>0$ \[ \chi(a_1+\epsilon h_1, \ldots,a_n + \epsilon h_n:
h_1,\ldots,h_n) \geq \mathbb P_{2 C \epsilon\sqrt{n}}(a_1,\ldots,a_n) + n 
\log \epsilon + \underline{\chi}(h_1,\ldots,h_n).\]
Dividing by $| \log \epsilon|$ on both sides, taking a $\limsup$ as
$\epsilon \rightarrow 0,$ and adding $n$ to both ends of the inequality
above yields 

\begin{eqnarray*} n + \limsup_{\epsilon \rightarrow 0} \frac
{\chi(a_1+\epsilon h_1, \ldots, a_n + \epsilon h_n: h_1,\ldots,
h_n)}{|\log \epsilon|} & \geq & \limsup_{\epsilon \rightarrow 0} \frac
{\mathbb P_{2 C \epsilon \sqrt{n}} (a_1,\ldots,a_n)}{|\log \epsilon|} \\ & 
= & \limsup_{\epsilon \rightarrow 0} \frac {\mathbb P_{\epsilon}
(a_1,\ldots,a_n)}{|\log \epsilon|}.  \end{eqnarray*}

\vspace{.2in}

	For the reverse inequality suppose $2 \leq m \in \mathbb N$ and 
$ \frac {1}{2(C+1)} > \epsilon > \sqrt{\gamma} > 0, R> \max_{1 \leq 
j \leq n} \{\|a_j\|\}.$  For each $k \in \mathbb N$ 
find an packing by open $\rho_k$ $\epsilon$-balls of 
$\Gamma_{R+1}(a_1, \ldots, a_n; 
m,k,\gamma)$ with maximum cardinality.  Denote the set of centers of the 
balls of this particular 
collection by $\Omega_k.$  Clearly
	
\begin{eqnarray*} \Gamma_{R+ \frac{1}{2}, \frac{1}{2}}\left 
(a_1 + \epsilon h_1, \ldots, a_n + \epsilon h_n: \epsilon h_1, 
\ldots, \epsilon h_n; m,k,\frac {\gamma}{2^m}\right) & \subset & 
\mathcal N_{2C\epsilon \sqrt{n}}(\Gamma_{R+1}(a_1, \ldots, a_n;m,k, 
\gamma)) 
\\ & \subset & \mathcal N_{3 C \epsilon \sqrt{n}}(\Omega_k) 
\end{eqnarray*}

\noindent where $\Gamma_{r + \frac{1}{2}, \frac{1}{2}}( \cdot )$ denotes 
the microstate space of $2n$-tuples such that the first $n$ entries have 
operator norms no larger than $r + \frac{1}{2}$ and the last $n$ entries 
have operator norms no larger than $\frac{1}{2}$ (see [4] for this 
technical modification).  $\mathcal N_{\epsilon}$ is taken with respect 
to the metric space $(M^{sa}_k(\mathbb C))^n$ with the $\rho_k$ metric.  
It follows that $\chi_{R + \frac{1}{2}, 
\frac{1}{2}}\left 
(a_1+ \epsilon h_1, \ldots, a_n + \epsilon h_n: \epsilon h_1, \ldots, 
\epsilon h_n; m, \frac {\gamma}{2^m} \right)$ is dominated by 

\begin{eqnarray*}\nonumber & & \limsup_{k \rightarrow \infty} 
\left [ k^{-2} \cdot \log (\text{vol}(\mathcal N_{3 C \epsilon 
\sqrt{n}}(\Omega_k))) 
+ \frac{n}{2} \cdot \log k \right] \\  & \leq & 
\limsup_{k \rightarrow \infty} \left [ k^{-2} \cdot \log 
\left(|\Omega_k| \cdot \frac {\pi^{\frac {nk^2}{2}} \cdot 
(3 C \epsilon \sqrt{nk})^{nk^2})}{\Gamma\left(\frac {nk^2}{2}+1 \right)} 
\right ) + \frac {n}{2} \cdot \log k \right ] \\  & \leq & 
\limsup_{k \rightarrow \infty} k^{-2} \cdot \log (|\Omega_k|) 
+ \limsup_{k \rightarrow \infty} \left [ n \log (3C \epsilon \sqrt{nk\pi}) 
- k^{-2} \cdot \log \left ( \frac {nk^2}{2e} \right)^{\frac {nk^2}{2}} + 
\frac {n}{2} \cdot \log k \right] \\ & = & \limsup_{k \rightarrow 
\infty} k^{-2} \cdot \log(|\Omega_k|) + \limsup_{k \rightarrow \infty} 
\left ( n\log (3 C \epsilon \sqrt{n\pi}) - n \log 
\left(\frac{k\sqrt{n}}{\sqrt{2e}} 
\right ) + n \log k \right) \\ & = & \limsup_{k \rightarrow \infty}
 k^{-2} \cdot \log(|\Omega_k|) + n \log (3C \epsilon \sqrt{2 \pi e})  \\
& = & \mathbb P_{\epsilon, R+1}\left (a_1, \ldots, a_n;m,
\frac{\gamma}{2^m} \right) + n \log(3C \epsilon \sqrt{2 \pi e}). 
\end{eqnarray*} 

\noindent This being true for any $2 \leq m \in \mathbb N,$
$\frac{1}{2(R+1)} > \epsilon > \sqrt {\gamma} > 0,$ and $R> \max_{1 \leq
j\leq n}\{\|a_j\|\}$ it follows that for sufficiently small $ \epsilon >0$
\begin{eqnarray*} \chi (a_1 + \epsilon h_1, \ldots, a_n + \epsilon h_n :
h_1, \ldots, h_n)  & = & \chi_{R+ \frac{1}{2}, \frac{1}{2}}(a_1 + \epsilon
h_1, \ldots, a_n + \epsilon h_n : \epsilon h_1 , \ldots, \epsilon h_n) \\
& \leq & \mathbb P_{\epsilon} (a_1, \ldots, a_n) + n \log \epsilon + n
\log(3C \sqrt{2 \pi e}).  \end{eqnarray*} Dividing by $|\log \epsilon|,$
taking a $\limsup$ as $\epsilon \rightarrow 0,$ and adding $n$ to both
sides of the inequality above yields \[ n + \limsup_{\epsilon \rightarrow
0} \frac {\chi(a_1 + \epsilon h_1, \ldots, a_n + \epsilon h_n : h_1,
\ldots, h_n)} {|\log \epsilon|} \leq \limsup_{\epsilon \rightarrow 0}
\frac {\mathbb P_{\epsilon}(a_1, \ldots, a_n)}{| \log \epsilon|}. \]
\end{proof}

\vspace{.15in} 

\begin{remark} Suppose $b_1,\ldots, b_p$ are contained in the
strongly closed algebra generated by the $a_i$ and $R >0$ is strictly
greater than the operator norm of any $a_i$ or $b_j.$ The proof shows that
the quanitity \[ \limsup_{\epsilon \rightarrow 0} \frac {\mathbb
K_{\epsilon, R}(a_1,\ldots, a_n:b_1,\ldots, b_p)}{ | \log \epsilon|} =
\limsup_{\epsilon \rightarrow 0} \frac {\mathbb P_{\epsilon,
R}(a_1,\ldots,a_n:b_1,\ldots, b_p)}{| \log \epsilon|}\]

\noindent equals \[ n + \limsup_{\epsilon \rightarrow 0}
\frac{\chi(a_1+\epsilon h_1, \ldots, a_n+\epsilon h_n: h_1,\ldots,
h_n)}{|\log \epsilon|}.\] \end{remark} 

\vspace{.15in}

Recall that by [3] and [5] if $\{s_1,\ldots,s_n\}$ is a free semicircular
family, then $\chi(s_1,\ldots,s_n) = \underline{\chi}(s_1,\ldots,s_n) > -
\infty.$ Thus we have by the lemma:

\begin{corollary}  
\[ \delta_0(a_1,\ldots,a_n) = 
\limsup_{\epsilon \rightarrow 0} 
\frac{\mathbb P_{\epsilon}(a_1,\ldots,a_n)}{|\log \epsilon|} 
= \limsup_{\epsilon \rightarrow 0} \frac {\mathbb 
K_{\epsilon}(a_1,\ldots,a_n)}{|\log \epsilon|}.\]
\end{corollary}

\hspace{.1in}

Both descriptions of $\delta_0,$ either in terms of volumes of
$\epsilon$-neighborhoods or in terms of packing numbers, can be useful.  
In the presence of freeness or in the situation with one random variable
it is fruitful to use the $\epsilon$-neighborhood description as
Voiculescu did ([3]).  On the other hand when computing $\delta_0$ in some
examples it is convenient to use the uniform packing description and this
was the implicit attitude taken towards $\delta_0$ in [2].  The packing
formulation also comes in handy when proving formulas for generators of
$M$ when $M$ has a simple algebraic decomposition into a tensor product of
a von Neumann algebra $N$ with the $k\times k$ matrices or into a direct
sum of algebras.

\end{document}